\documentclass[11pt]{amsart}
\usepackage{amsmath, amsthm, latexsym, amssymb, amsfonts,
epsfig}

\newenvironment{pf}{\proof[\proofname]}{\endproof}
\theoremstyle{plain}
\newtheorem{Th}{Theorem}[section]

\newtheorem{Lemma}[Th]{Lemma}
\numberwithin{equation}{section} 
\numberwithin{figure}{section} 
\theoremstyle{definition}

\newtheorem{Def}[Th]{Definition}

\newcommand{\res}{\text{res}}
\newcommand{\cal}[1]{\mathcal{#1}}
\newcommand{\C}{\mathbb C}
\newcommand{\Z}{\mathbb Z}
\newcommand{\R}{\mathbb R}

\renewcommand{\Re}{\text{\rm Re}}

\newcommand{\w}{\widetilde}

\newcommand{\rl}[1]{Lemma~\ref{L:#1}}

\newcommand{\rf}[1]{Figure~\ref{F:#1}}

\newcommand{\rt}[1] {Theorem~\ref{T:#1}}

\begin{document}


\title{Exponential Gelfond-Khovanskii formula in dimension one}
\author{Evgenia Soprunova}
\address{Department of Mathematics and Statistics,
University of Massachusetts, Amherst, MA 01003}
\email{soproun@math.umass.edu}

\date{}
\begin{abstract}
Gelfond and Khovanskii found a formula for the sum of the values
of a Laurent polynomial over the zeros of a system of $n$ Laurent
polynomials in $(\C\setminus 0)^n$. We expect that a similar
formula holds in the case of exponential sums with real
frequencies. Here we prove such a formula in dimension 1.
\end{abstract}
\maketitle{}

\section{Introduction}
O.~Gelfond and A.~Khovanskii found a formula for the sum
of the values of a Laurent polynomial $g$ over the zeros of a system
$$f_1(z)=\dots=f_n(z)=0,\quad z\in(\C\setminus 0)^n,$$ where
$f_i$ are Laurent polynomials whose Newton polyhedra have generic
relative positions \cite{GKh1,GKh2}. After an exponential
change of variables we obtain a similar formula where $g$ and the $f_i$ are 
now exponential sums with rational frequencies. There is a lot of evidence that this formula
continues to hold if the frequencies are real. For example, when $g=1$ the formula follows
from combining two results: Gelfond's generalization of
Bernstein's theorem \cite{G2} and the new formula for mixed volume
\cite{Kh1}. In \cite{So} the formula is proved in the
case when the frequencies of the exponential sum $g$ are not commensurate with the 
frequencies of the system.

Here we prove the generalization of the Gelfond--Khovanskii
formula to exponential sums with real frequencies in dimension
one. Our argument is elementary.

A similar result is obtained by J. Ritt in \cite{Ritt}, where he
computes the average sum of the real parts of zeros of an
exponential sum. Ritt's result can be considered as an exponential
generalization of the Vieta formula for the product of zeros of a
polynomial.

This paper is a part of the author's Ph.D. thesis. I would like to
thank my thesis advisor Askold Khovanskii for stating the
problem and for his constant attention to this work.

\section{Algebraic case}
Recall that a Laurent polynomial $f$ in one complex variable $z$ is a finite
linear combination of monomials with integer exponents:
$$f(z)=\sum_k c_kz^k,\quad k\in\Z,\quad c_k\in\C.$$
We explain how to compute explicitly the sum of the values of a Laurent
polynomial $g$ over the zeros of a Laurent polynomial $f$ in
$\C\setminus\{0\}$, counting multiplicities. Our argument in the exponential case is a
generalization of the argument that we present here.

Consider a differential form $\omega=g\,df/f$. This form has simple poles
at the zeros $z_i$ of $f$ with the residue equal to $\mu_i g(z_i)$ where
$\mu_i$ is the multiplicity of the root $z_i$. Since the sum 
of the residues of $\omega$ is zero, the sum of the values of $g$ over the zeros
of $f$ in $\C\setminus\{0\}$ is equal to $-\res_{\infty}\,\omega-\res_{0}\,\omega.$

To compute these residues, we expand $gf'/f$ into Laurent series as follows. Let
$$f(z)=c_1z^{k_1}+\cdots+c_nz^{k_n},\quad k_1<\cdots<k_n,\quad k_i\in\Z.$$ Put $\w
f=f/c_1 z^{k_1}$. Since $\w f-1$ contains only positive powers 
of $z$,  we have $|\w f-1|<1$ in a neighborhood of zero.
Therefore, we
can expand  $1/\w f$ into series converging uniformly in this
neighborhood:
$$
1/{\w f}=1+(1-\w f)+(1-\w f)^2+\cdots.
$$
Since each power of $z$ appears in this series finitely many
times, it is a well-defined series in $z$. Multiplying this series by
$gf'/c_{1}z^{k_1}$ we obtain a series for $gf'/f$:
$$\frac{gf'}{f}
=\frac{gf'}{c_1z^{k_1}}(1+(1-\w f)+(1-\w f)^2+\cdots).
$$
 Therefore, the residue at zero
is equal to the coefficient $A_1$ of $1/z$ in this series.

Similarly, denote by $A_n$ the coefficient of
$1/z$ in the series
$$\frac{gf'}{c_nz^{k_n}}(1+(1-\w
f)+(1-\w f)^2+\cdots),$$ where this time $\w f=f/c_nz^{k_n}$.
Then the residue at infinity  is equal to $-A_n$.

Therefore, the sum of the values of $g$ over the zeros of
$f$ in $\C\setminus\{0\}$ is equal to $A_n-A_1$.

\section{Statement of result}
We are going to deal with exponential sums with real frequencies,
i.e. functions of the form
$$f(z)=c_1\exp 2\pi\alpha_1 z+\cdots +c_n\exp 2\pi\alpha_n z,$$
where the coefficients $c_i$ are nonzero complex numbers and the
frequencies $\alpha_i$ are increasing real numbers:
$\alpha_1<\dots<\alpha_n.$  If we factor out $c_1\exp 2\pi\alpha_1 z$
from $f$, the rest will be close to infinity as $\Re\,z\to\infty$
and close to 1 as $\Re\,z\to -\infty$. Therefore, all the zeros of $f$ lie in the vertical
strip $|\Re\,z|<B$ for some~$B$, where the choice of~$B$
depends only on the frequencies $\alpha _i$ and the absolute
values of the coefficients~$c_i$.

Let $g$ be an exponential sum with real frequencies. We
want to add up the values of $g$ over the zeros of
$f$. Since the number of zeros is infinite,
the sum over the zeros must be replaced with the result of
averaging $g$ over the zeros of $f$ along the imaginary
axis. More precisely, let $S(R)$ be the sum of the values of $g$ over those
zeros of $f$ (counting multiplicities) whose absolute value
of the imaginary part is less than $R$. Since all the zeros of $f$
belong to the vertical strip $|\text{Re}\,z|<B$, this sum is finite.

\begin{Def}
The {\it mean value} of $g$ over the zeros of $f$  
is the limit
$$M=\lim_{R\to\infty}\frac{S(R)}{2R}.$$
If $g=1$, $M$ is the {\it mean number} of zeros of $f$. 
\end{Def}

Let $\w f=f/c_k\exp 2\pi\alpha_k z$, where $k$ is either 1 or $n$.
The constant term of the exponential sum $\w f$
is equal to one, thus we can define the exponential series for $1/\w f$ by
the formula 
$$1/\w f=1+(1-\w f)+(1-\w f)^2+\cdots.$$ 
Since each exponent $\exp 2\pi \alpha z$ appears with a nonzero coefficient in a
finite number of terms, the coefficients of these series are
well-defined. Let $A_k$ be the constant term in the formal product
of this series and $gf'/c_k\exp(2\pi\alpha_k z)$ for $k=1,n$.

The following theorem is the generalization of the
Gelfond--Khovanskii formula to the case of exponential sums in
dimension one.

\begin{Th}\label{T:onedim}
Let $f$ be an exponential sum with real frequencies:
$$
f(z)=c_1\exp 2\pi\alpha_1 z+\cdots +c_n\exp 2\pi\alpha_n z, \quad \alpha_1<\dots<\alpha_n.
$$
Let $g$ be another exponential sum with real frequencies.
Then the mean value $M$ of $g$ over the zeros of $f$ is equal to
$$\frac{A_n-A_1}{2\pi},$$ where $A_n$ and $A_1$ are defined above.
In particular, the mean number of zeros of $f$, is equal to
$\alpha_n-\alpha_1$.

Furthermore, the mean value of $g=\exp
2\pi\alpha z$ over the zeros of $f$ can only be nonzero if
$\alpha$ belongs either to the non-positive semigroup $\cal A$
generated by the set $\{\alpha_1-\alpha_i\,|\,i=1\dots n\}$, or to the
non-negative semigroup $\cal B$ generated by the set
$\{\alpha_n-\alpha_i\,|\,i=1\dots n\}$. 
\end{Th}

\section{The plan of the proof}
We can assume that $f$ does not have zeros on $|\text{Im}\,z|=R$. We have
$$M=\lim_{R\to\infty}\frac{S(R)}{2R}=\lim_{R\to\infty}\frac{1}{4\pi iR}\int_{\Gamma_R}\frac{gf'}{f}dz,$$
where the integration is performed in the positive sense over the
rectangle $\Gamma_R$ bounded by the lines $\text{Re}\,z=B$, $\text{Im}\,z=R$, $\text{Re}\,z=-B$,
and $\text{Im}\,z=-R$. Denote the corresponding sides of $\Gamma_R$ by
$\gamma_+$, $\gamma_R$, $\gamma_-$, and $\gamma_{-R}$ (see \rf{contour}).
\setlength{\unitlength}{1mm}
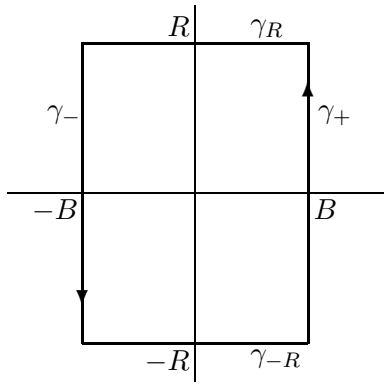
\begin{figure}
\begin{picture}(50,50)(-25,-25)
\thicklines
\put(-15,-20){\line(1,0){30}}
\put(15,-20){\line(0,1){40}}
\put(-15,-20){\line(0,1){40}}
\put(-15,20){\line(1,0){30}}
\put(15,10){\vector(0,1){5}}
\put(-15,-10){\vector(0,-1){5}}
\thinlines
\put(-25,0){\line(1,0){50}}
\put(0,-25){\line(0,1){50}}
\put(15.4,-3.4){$B$}
\put(-22,-3.4){$-B$}
\put(-4,21){$R$}
\put(-7,-23.5){$-R$}
\put(7,21.4){$\gamma_R$}
\put(7,-22.5){$\gamma_{-R}$}
\put(16,10){$\gamma_+$}
\put(-20,10){$\gamma_-$}
\end{picture}
\caption{Contour $\Gamma_R$}\label{F:contour}
\end{figure}

The proof of \rt{onedim} will
be divided in two parts between the following two sections as
follows:
\begin{enumerate}
\item We show that the integral over the horizontal part of $\Gamma_R$ is
bounded as $R\to\infty$, and consequently only the integrals over
the vertical segments contribute to the mean value.
\item We compute the contribution to
the mean value of the integrals over the vertical segments.
\end{enumerate}

\section{The integral over the horizontal part of $\Gamma_R$ is bounded}
For each zero of the exponential sum 
$f(z)=c_1\exp 2\pi\alpha_1z+\cdots+c_n\exp 2\pi\alpha_n z$ 
we delete from $|\text{Re}\,z|\leq B$ 
an open disk of radius $r>0$ centered at that zero, and
call the resulting region $D_f$. In the lemma below we show that
$f$ is separated from~0 on~$D_f$.

\begin{Lemma}\label{L:separ}
For $D_f$ defined above there exists $\delta>0$ such that for all $z\in D_f$
$$|f(z)|\geq \delta>0.$$
\end{Lemma}
\begin{pf}
Let $\cal F$ be the family of exponential sums that have the same
frequencies $\alpha_i$ as~$f$, and each of its coefficients
$d_i$ has the same absolute value as the corresponding $c_i$ in~$f$. 
That is, each $h\in\mathcal F$ is of the form
$$h(z)=d_1\exp 2\pi\alpha_1z+\dots +d_n\exp 2\pi \alpha_n z,$$
where $|d_i|=|c_i|.$ This family is invariant under translations
along the imaginary axis. Due to our choice of $B$, all the zeros of
each of the functions in $\cal F$ belong to $|\text{Re}\,z|< B$.

For each element $h$ of $\cal F$ we define $D_h$ to be the strip
$|\text{Re}\,z|\leq B$ with deleted open disks of radius $r$
centered at the zeros of $h$. Let
$$\delta_h=\min_{z\in D_h}|h(z)|,\quad \delta=\inf_{h\in\cal F}\delta_h.$$
We need to show that $\delta>0$.

Assume $\delta=0$. It means that there exist a
sequence of functions $h_l\in\cal F$ and a sequence of points
$z_l\in D_{h_l}$ such that $h_l(z_l)\to 0$, as $l\to\infty$.
Since $\cal F$ is invariant under translations along the imaginary
axis, we can assume that all $z_l$ belong to the compact set
$$V=\{|\text{Re}\,z|\leq B,\, |\text{Im}\,z|\leq 1\}.$$
Each $h\in\cal F$ can be identified with
the collection of its coefficients $(d_1,\dots,d_n)\in\C^n$. Then
the convergence in $\C^n$ corresponds to the uniform convergence
on $V$. The family $\cal F$ is a product of $n$
circles and, therefore, is a compact set.

Consider the sequence $(h_l,z_l)$ in $\cal F\times V$.
Since  $\cal F\times V$ is compact, we can assume (choosing a subsequence 
if necessary) that the sequence  $(h_l, z_l)$ converges
to some $(h,z_0)\in\cal F\times V$. This means that $z_l$ converges to $z_0$, 
and $h_l$ converges to $h$ uniformly on $V$. Since
$h_l(z_l)\to 0$, we have $h(z_0)=0$.

By the Hurwitz theorem, each function $h_l$ with $l$ big
enough has a zero in the disk $\{|z-z_0|<\rho\}\subset\{ |\text{Re}\,z|\leq B\}$ 
with $\rho<r/2$. Therefore, we can find $l$ such
that $z_l$ belongs to the disk and $h_l$ has a zero in
this disk. But then the distance between $z_l$ and this zero
of $h_l$ is less than $r$, which contradicts $z_l\in
D_{h_l}$.

This contradiction proves that $\delta>0$.
\end{pf}

We will need the following result from the theory of
fewnomials, which is a consequence of the complex Rolle theorem.

\begin{Th}\label{T:fewnom}\cite{Ya}
Let $f$ be an exponential sum of the form
$$f(z)=c_1\exp 2\pi\alpha_1z+\dots +c_n\exp 2\pi\alpha_n z,$$
where the coefficients $c_i$ are complex numbers and the
frequencies $\alpha_i$ are increasing real numbers:
$\alpha_1<\dots<\alpha_n.$ Then $f$ has less than $n$ zeros in a
horizontal strip of width less than $1/(\alpha_n-\alpha_1).$
\end{Th}

Recall that $\gamma_R\cup\gamma_{-R}$ is the horizontal part of the
contour $\Gamma_R$.

\begin{Lemma}
The absolute value of the integral
$$\int_{\gamma_R}\frac{gf'}{f}dz$$
is bounded for all $R$ such that there are no zeros of $f$ on $\gamma_R$.
\end{Lemma}

\begin{pf}
Choose $r< 1/4n(\alpha_n-\alpha_1)$. For each
zero of $f$ delete an open disk of radius $r$ centered at
that zero. The projection of these deleted disks to the imaginary
axis does not entirely cover the vertical segment
$$\Big\{\text{Re}\,z=0,\,\big|\,\text{Im}\,z-R\,\big|<\frac{1}{4(\alpha_n-\alpha_1)}\Big\}.
$$
Indeed, the only disks whose projections could touch this segment
are the ones whose center $z$ satisfies
$$\big|\,\text{Im}\,z-R\,\big|<\frac{1}{2\,(\alpha_n-\alpha_1)}.$$ 
By \rt{fewnom}, there are
less than $n$ such disks. Hence the length of the projection
is less than $2rn<1/2\,(\alpha_n-\alpha_1)$, which is
equal to the length of the vertical segment. Therefore, we can choose $R'$ on this segment so that
$\gamma_{R'}$ does not meet the deleted disks. Since by \rl{separ}
$f$ is separated from zero on $\gamma_{R'}$, the integral over
$R'$ is bounded.

It remains to show that the difference of the integrals over
$\gamma_R$ and $\gamma_{R'}$ is bounded. Consider the contour
$\gamma$ formed by $\gamma_R$, $\gamma_{R'}$, and the vertical
segments connecting their endpoints. The integral over $\gamma$ is
equal to the sum of the values of $g$ over the zeros of $f$
that lie inside $\gamma$, times $2\pi i$. This sum is bounded since $f$ has less than
$n$ zeros inside, and $g$ is bounded on $|\text{Re}\,z|<B$.
Choosing $B$ so big that all the zeros of $f$ belong to
$|\text{Re}\,z|<B/2$, and requiring that $r<B/2$, we ensure that
the deleted disks do not intersect the vertical part of $\gamma$, and the
integral over the vertical segments is bounded. We
proved that the integral
$$\int_{\gamma_R}\frac{gf'}{f}dz$$
is bounded.
\end{pf}

This lemma implies that
$$\lim_{R\to\infty}\frac{1}{4\pi iR}\int_{\gamma_R}\frac{gf'}{f}dz=0.
$$
Therefore, only the integrals over the vertical part of $\Gamma_R$ can
contribute to the mean value.

\section{The integral over the vertical part of $\Gamma_R$.  }

Here we  explain how to compute the limit
$$\lim_{R\to\infty}\frac{1}{4\pi iR}\int\frac{gf'}{f}dz,
$$
where the integration is performed over the vertical part of
$\Gamma_R$. 

We deal with $\gamma_-$ first. Recall that
$f(z)=c_1\exp2\pi\alpha_1z+ \dots+c_n\exp2\pi\alpha_nz$, where
$\alpha_1<\dots<\alpha_n.$ Let $\w f=f/c_1\exp2\pi\alpha_1
z$. Then the constant term of $\w f$ is equal to 1. On
$\text{Re}\,z=-B$,
$$|\w f-1|=\left|\sum_{j=2}^n \frac{c_j}{c_1}\exp2\pi(\alpha_j-\alpha_1)z\right|
\leq
\sum_{j=2}^n\frac{|c_j|}{|c_1|}\exp(-2B\pi(\alpha_j-\alpha_1)),$$
which is less than 1 if $B$ is big enough since
$\alpha_j>\alpha_1$. Therefore, we can expand $gf'/f$ into
geometric series converging uniformly on $\text{Re}\,z=-B$:
$$
g\frac{f'}{f}=
g\frac{f'}{c_1\exp2\pi\alpha_1 z}\left(1+(1-\w f)+(1-\w
f)^2+\dots\right),
$$
which is a series in exponents, since every exponent appears only
finitely many times.  We integrate this series over $\gamma_-$, divide by
$4\pi i R$ and find the limit as $R\to\infty$ term-by-term.

Notice that only the constant term in this series can give a
nonzero limit. If this constant term is $A_1$, then the
corresponding limit is $-A_1/2\pi$. If $g(z)=\exp2\pi\alpha z$,
then the only values of $\alpha$ for which there is a nonzero
constant term in the series are those that belong to the non-positive semigroup
$\cal A$ generated by $\{\alpha_1-\alpha_i\,|\,i=1,\dots,n\}$.

It remains to repeat this procedure for $\gamma_+.$ Let $\w
f=f/c_n\exp2\pi\alpha_n z$, and expand $gf'/f$ into 
exponential series converging uniformly on $\text{Re}\, z=B$:
$$\frac{gf'}{f}=
\frac{gf'}{c_n\exp2\pi\alpha_n z}\left(1+(1-\w f)+(1-\w
f)^2+\dots\right).
$$
If the constant term of this series is $A_n$, then the
corresponding limit is $A_n/2\pi$. If $g(z)=\exp2\pi\alpha
z$, then the only values of $\alpha$ for which there is a nonzero
constant term in the series are those that belong to the non-negative
semigroup
$\cal B$ generated by $\{\alpha_n-\alpha_i\,|\,i=1,\dots,n\}$.

We proved that the mean value $M$ of $g$ over the zeros of
$f$ is equal to $$\frac{A_n-A_1}{2\pi}.$$ If $g=1$, it is easy to see 
that $A_n=2\pi\alpha_n$ and $A_1=2\pi\alpha_1$. Hence the mean number
of zeros of $f$ is equal to $\alpha_n-\alpha_1$.


\begin{thebibliography}{99}

\bibitem[1]{GKh1} O. Gelfond, A. Khovanskii,
{\em Newton polyhedra and Grothendieck residues} (in Russian),
Dokl. Akad. Nauk, {\bf 350}, no. 3 (1996), 298--300.

\bibitem[2]{GKh2} O. Gelfond, A. Khovanskii,
{\em Toric geometry and Grothendieck  residues}, Moscow Mathematical Journal, Vol. {\bf 2},
 no. 1 (2002), 99--112.

\bibitem[3]{G1} O. Gelfond, {\em Zeros of systems of quasiperiodic
polynomials}, FIAN preprint, No. 200 (1978).

\bibitem[4]{G2} O. Gelfond,
{\em The mean number of roots of systems of holomorphic almost
periodic equations.} (Russian) Uspekhi Mat. Nauk {\bf 39} (1984), no.
1(235), 123--124.

\bibitem[5]{Kh1} A. Khovanskii,
{\em Newton polyhedra, a new formula for mixed volume, product of
roots of a system of equations}, Fields Inst. Comm., Vol.{\bf 24}
(1999), 325--364.

\bibitem[6]{Ya} A. Khovanskii, S. Yakovenko,
{\em Generalized Rolle Theorem in $\R^n$ and $\C$}, J. Dynam.
Control Systems {\bf 2} (1996), no. 1, 103--123.

\bibitem[7]{Ritt} J. Ritt,
{\em On the zeros of exponential polynomials}, Transactions of The
American Mathematical Society, Volume 31, Issue 4 (Oct., 1929),
680-686.

\bibitem[8]{So} E. Soprunova,
{\em Zeros of systems of exponential sums and trigonometric
polynomials}, preprint, 2002.


\end{thebibliography}
\end{document}